\documentclass[12pt]{amsart}
\usepackage{amsfonts,amssymb,color}
\usepackage[mathscr]{eucal}
\usepackage{amsmath, amsthm}
\usepackage{mathrsfs}
\usepackage{amsbsy}
\usepackage{dsfont}
\usepackage{bbm}
\makeindex
\input xypic
\xyoption{all}

\begin{document}
\baselineskip = 5mm
\newcommand \lra {\longrightarrow}
\newcommand \hra {\hookrightarrow}
\newcommand \ZZ {{\mathbb Z}} 
\newcommand \NN {{\mathbb N}} 
\newcommand \QQ {{\mathbb Q}} 
\newcommand \RR {{\mathbb R}} 
\newcommand \CC {{\mathbb C}} 
\newcommand \bcA {{\mathscr A}}
\newcommand \bcB {{\mathscr B}}
\newcommand \bcC {{\mathscr C}}
\newcommand \bcD {{\mathscr D}}
\newcommand \bcE {{\mathscr E}}
\newcommand \bcF {{\mathscr F}}
\newcommand \bcT {{\mathscr T}}
\newcommand \bcU {{\mathscr U}}
\newcommand \bcX {{\mathscr X}}
\newcommand \bcY {{\mathscr Y}}
\newcommand \bcZ {{\mathscr Z}}
\newcommand \C {{\mathscr C}}
\newcommand \im {{\rm im}}
\newcommand \Hom {{\rm Hom}}
\newcommand \colim {{{\rm colim}\, }} 
\newcommand \End {{\rm {End}}}
\newcommand \coker {{\rm {coker}}}
\newcommand \id {{\rm {id}}}
\newcommand \supp {{\rm {Supp}}\, }
\newcommand \CHM {{\mathscr C\! \mathscr M}}
\newcommand \DM {{\mathscr D\! \mathscr M}}
\newcommand \MM {{\mathscr M\! \mathscr M}}
\newcommand \uno {{\mathbbm 1}}
\newcommand \Le {{\mathbbm L}}
\newcommand \PR {{\mathbb P}} 
\newcommand \AF {{\mathbb A}} 
\newcommand \Spec {{\rm {Spec}}}
\newcommand \Pic {{\rm {Pic}}}
\newcommand \Alb {{\rm {Alb}}}
\newcommand \Corr {{Corr}}
\newcommand \Sym {{\rm {Sym}}}
\newcommand \cha {{\rm {char}}}
\newcommand \tr {{\rm {tr}}} 
\newcommand \gm {{\mathfrak {m}}}
\newcommand \gp {{\mathfrak {p}}}
\def\blue {\color{blue}}
\def\red{\color{red}}
\newtheorem{theorem}{Theorem}
\newtheorem{lemma}[theorem]{Lemma}
\newtheorem{sublemma}[theorem]{Sublemma}
\newtheorem{corollary}[theorem]{Corollary}
\newtheorem{example}[theorem]{Example}
\newtheorem{exercise}[theorem]{Exersize}
\newtheorem{proposition}[theorem]{Proposition}
\newtheorem{remark}[theorem]{Remark}
\newtheorem{notation}[theorem]{Notation}
\newtheorem{definition}[theorem]{Definition}
\newtheorem{conjecture}[theorem]{Conjecture}
\newtheorem{claim}[theorem]{Claim}
\newenvironment{pf}{\par\noindent{\em Proof}.}{\hfill\framebox(6,6)
\par\medskip}
\title[Motives and cycles on threefolds over a field]
{\bf Motives and representability of algebraic cycles on threefolds over a field}
\author{S. Gorchinskiy, V. Guletski\u \i }

\date{12 November 2012}

\begin{abstract}
We study algebraic cycles on threefolds and finite-dimen\-sionality of their motives with coefficients in $\mathbb Q$. We decompose the motive of a non-singular projective threefold $X$ with representable algebraic part of $CH_0(X)$ into Lefschetz motives and the Picard motive of a certain abelian variety, isogenous to the Griffiths' intermediate Jacobian $J^2(X)$ when the ground field is $\mathbb C$. In particular, it implies motivic finite-dimensionality of Fano threefolds over a field. We also prove representability of zero-cycles on several classes of threefolds fibred by surfaces with algebraic $H^2$. This gives another new examples of three-dimensional varieties whose motives are finite-dimensional.
\end{abstract}

\maketitle

\section{Introduction}
\label{s-intro}

It is well known that algebraic cycles in codimension greater than
one are hard to understand even for surfaces. A general
Bloch-Beilinson's philosophy says that the Chow groups of all
smooth projective varieties should be filtered, and the graded
pieces of those filtrations should be isomorphic to Ext-groups in a
certain conjectural category of mixed motives over the ground field.
Finite-dimensional motives, which had been introduced by S.Kimura in
\cite{Kimura}, shed a new light on the motivic picture of
intersection theory. If $X$ is a surface over an algebraically
closed field with algebraic second cohomology group, the Chow group
of zero cycles on $X$ is representable if and only if the
motive $M(X)$ is finite-dimensional, see \cite{GP2}, Theorem 7. The
purpose of the present paper is to investigate some interesting
links between motivic finite-dimensionality and representability of algebraic cycles on three-dimensional varieties over an arbitrary algebraically closed ground field $k$. In particular, we will show that the motive of a smooth projective threefold $X$ over $k$ can be decomposed into a sum of Lefschetz and Abelian motives if and only if zero-cycles are representable on $X$ (Theorem \ref{t-zerocycles}). To some extent, this result can be considered as a most possible analog of Theorem 7 from \cite{GP2}. Notice that representability of zero-cycles on Fano threefolds was proved uniformly by Koll\'ar in \cite{Kollar}, which was a generalization of the results of Bloch and Murre, see a survey in \cite{Mu3}. Joint with our result it gives motivic finite-dimensionality for all Fano threefolds. Then we prove representability of algebraically trivial algebraic cycles on several new types of threefolds, in characteristic zero, fibred by surfaces with algebraic $H^2$ (Theorem \ref{main2}). After all, this gives another series of threefolds whose motives are finite-dimensional over a field of characteristic zero (Theorem \ref{findimconcrete} and Theorem \ref{findimconcrete2}). One can hope that such fund of threefolds with known finite-dimensionality of their motives can be useful in approaching the problem of motivic finite-dimensionality for $K3$-surfaces.

The paper is organized as follows. First we recall Chow motives over a base and prove a sort of ``Manin's lemma" which is interesting itself and needed in further computations. Then we recall some necessary things about weak representability of algebraic cycles, and then develop a motivic vision of that effect suitable in our approach. All of this covers Sections 2 - 4. In Section 5 we state and prove Theorem \ref{t-zerocycles} bringing a motivic criteria for representability of algebraic cycles on threefolds over an arbitrary algebraically closed field. In Section \ref{s-fibred} we investigate threefolds fibred by surfaces with $p_g=0$ and show when algebraic cycles on such threefolds are representable, so their motives are finite-dimensional. Finally, in the last section we apply the previous results in the very concrete setting and prove representability of algebraic cycles and motivic finite-dimensionality for threefolds fibred by Enriques and hyperelliptic surfaces.

\section{Some motivic lemma}
\label{gerstenization}

Throughout this paper we use various categories of motives over different base schemes and with coefficients in $\QQ$ and $\ZZ$.

Let $k$ be an algebraically closed field, and let $X$ be a smooth
variety over $k$. As we are going to work with motivic finite-dimensionality, all Chow groups $CH^i(X)$ will be mostly with coefficients in $\QQ $, and we will use the notation $CH_{\ZZ }^i(X)$ for Chow groups with integral coefficients. For any codimension $i$ let $A^i(X)$ be a subgroup in $CH^i(X)$ generated by cycles algebraically equivalent to zero. Respectively, $A^i_{\ZZ }(X)$ is a subgroup of algebraically trivial cycle
classes in $CH^i_{\ZZ }(X)$.

Let $\CHM (S)$ be the category of relative Chow motives over an arbitrary irreducible smooth variety $S$ over $k$ constructed either contravariantly or covariantly. Here we recall the construction of $\CHM (S)$, for more information the reader may consult the papers \cite{DM} and \cite{CH}. In contravariant notation, if $X$ and $Y$ are two varieties over $S$ such that the structure morphisms $X\to S$ and $Y\to S$ are smooth and projective, and $X=\cup _jX_j$ are the connected components
of $X$, then
$$
\Corr _S^m(X,Y)=\oplus _jCH^{e_j+m}(X_j\times _SY)
$$
is a group of relative correspondences of degree $m$ from $X$ to $Y$
over $S$, where $e_j$ is the relative dimension of $X_j$ over $S$.
Given a morphism $f:X\to Y$ over $S$, the transpose $\Gamma _f^t$ of
its graph $\Gamma _f$ is in $\Corr _S^0(X,Y)$. For any two
correspondences $f:X\to Y$ and $g:Y\to Z$ their composition $g\circ
f$ is defined by a standard formula
$$
g\circ f={p_{13}}_*(p_{12}^*(f)\cdot p_{23}^*(g))\; ,
$$
where the central dot denotes the intersection of
cycle classes in Fulton's sense, \cite{Ful}, and the projections are
projections of a fibred product over $S$. Objects in $\CHM (S)$ are
triples $(X/S,p,n)$ where $p\in \Corr _S^0(X,X)$ is a projector
(i.e., idempotent) and $n$ is an integer. For motives $M=(X/S,p,m)$ and
$N=(Y/S,q,n)$, we have
$$
\Hom(M,N)=q\circ\Corr_S^{n-m}(X,Y)\circ p\; .
$$
For any $X/S$ its motive $M(X/S)$ is defined by the relative
diagonal $\Delta _{X/S}$, $M(X/S)=(X/S,\Delta _{X/S},0)$, and for
any morphism $f:X\to Y$ over $S$ the correspondence $\Gamma _f^{t}$
defines a morphism $M(f):M(Y/S)\to M(X/S)$. Thus, we have a
contravariant functor from the category of smooth projective
varieties over $S$ to $\CHM(S)$.

It is a substantial matter that the category $\CHM (S)$ is rigid
with a tensor product satisfying the formula
  $$
  (X/S,p,m)\otimes (Y/S,q,n)=(X\times _SY,p\otimes_S q,m+n)\; .
  $$
The scheme $S/S$ indexed by $0$ gives the unite motive $\uno $, and being indexed by $-1$ it gives the Lefschetz motive $\Le =(S,\Delta_{S/S},-1)$. Later on, for short, we will write $\Le ^n$ instead of tensor powers $\Le ^{\otimes n}$. The duality in the category $\CHM(S)$ is defined as follows:
$$
(X/S,p,m)^{\vee}=(X/S,p^{\rm t},d-m)\; ,
$$
where $X/S$ is of pure relative dimension $d$ over $S$.

Notice that any degree zero correspondence $f:X\to Y$ between smooth projective varieties over $S$ acts on (absolute) Chow groups
    $$
    f_*:CH^i(X)\lra CH^i(Y)
    $$
by the formula
  $$
  f_*(z)=(p_2)_*(p_1^*(z)\cdot f)\; ,
  $$
for any $z\in CH^i(X)$. For a motive $M=(X,p,m)$ in $\CHM(S)$ its Chow groups are defined as follows:
$$
CH^i(M)=\im (p_*:CH^{i+m}(X)\lra CH^{i+m}(X)).
$$
Analogously, groups of algebraically trivial cycles of $M$ are defined by the formula
$$
A^i(M)=\im (p_*:A^{i+m}(X)\lra A^{i+m}(X))\; .
$$

Fix a prime $l\ne\cha(k)$ and let
  $$
  H^i(X)(j)=H^i_{\acute e t}(X,\QQ_l(j))
  $$
be the $l$-adic \'etale cohomology group for variety $X$ over the field $k$. Correspondences act also on cohomology, and cohomology groups of a motive $M=(X,p,m)$ in $\CHM(S)$ are defined by the formula:
$$
H^i(M)(j)=p_*H^{i+2m}(X)(j+m)\; .
$$

If $f:S\to S'$ is a morphism of base schemes over $k$, then $f$ gives a base change tensor functor $f^*:\CHM (S')\lra \CHM (S)$. If $F$ is a field, we will write $\CHM (F)$ for the category of Chow motives over $\Spec (F)$.

Sometimes we will use the category of integral Chow motives
$\CHM_{\ZZ}(S)$. These are defined in the same way as rational ones but all Chow groups must be taken with integral coefficients. For a variety $X/S$, let $M_{\ZZ}(X/S)$ be the integral motive $(X/S,\Delta_{X/S},0)$ in the
category $\CHM_{\ZZ}(S)$.

By Manin's identity principle, a motive $M$ in $\CHM (k)$ is trivial if and only if $CH^*(M\otimes M(X))=0$ for any smooth projective variety $X$ over $k$. We will need a slightly different assertion:

\begin{lemma}
\label{l-gersten} Let $\Omega $ be a universal domain over $k$, i.e.
an algebraically closed field of infinite transcendence degree over
$k$. Let $M\in \CHM (k)$ and assume that $CH^*(M_{\Omega })=0$,
where $M_{\Omega }$ is the pullback of $M$ to $\CHM (\Omega )$
induced by the morphism $\Spec (\Omega )\to \Spec (k)$. Then it
follows that $M=0$ in $\CHM (k)$.
\end{lemma}

\begin{proof}
Let $M=(V,q,n)$ be an object in $\CHM (k)$ and assume that the group $CH^*(M_{\Omega })$ is trivial in $\CHM (\Omega )$. By Manin's principle, in order to prove that $M=0$ in $\CHM (k)$ it is enough to show that $CH^*(M\otimes M(X))=0$ for any smooth projective and irreducible variety $X$ over $k$. But in fact we will prove more.

Let $X$ be any irreducible variety over $k$. The projections
$p_{13}:V\times V\times X\to V\times X$ and $p_{23}:V\times V\times
X\to V\times X$ are flat and proper morphisms of algebraic
varieties. Let us define an action
$$
q\otimes X: CH^i(V\times X)\to CH^i(V\times X)
$$
by the formula:
$$
(q\otimes X)(a)={p_{23}}_*(p_{13}^*(a)\cdot_{p_{12}}q)\; .
$$
Here we use an operation $*\cdot_{p_{12}}*$ defined in
\cite[Ch.8]{Ful}, i.e. intersection with respect to the projection $p_{12}$. This is possible because $p_{12}$ is a projection onto a smooth variety. Notice that, as $q$ is an idempotent, so is $q\otimes X$.

Let now $U$ be a Zariski open subset in $X$ and let $Z=X-U$. The localization exact sequence
$$
CH_{e+d-i}(V\times Z)\lra CH^i(V\times X)\lra CH^i(V\times U)\to 0\; ,
$$
where $e=\dim (V)$ and $d=\dim (X)$, gives an exact sequence
$$
\oplus _UCH_{e+d-i}(V\times Z)\lra CH^i(V\times X)\lra \colim CH^i(V\times U)\to 0\; .
$$
Since the colimit is, actually, isomorphic to $CH^i(V_K)$, where $K=k(X)$, we get an exact sequence
$$
\oplus _UCH_{e+d-i}(V\times Z)\lra CH^i(V\times X)\lra CH^i(V_K)\to 0\; .
$$
Applying the projectors $q\otimes Z$, $q\otimes X$ and $q\otimes K$, using their compatibility with pull-backs and push-forwards and passing to images of idempotents we obtain an exact sequence
$$
\oplus _U\im (q\otimes Z)\lra \im (q\otimes X)\lra CH^i(M_K)\to 0\; .
$$

Now we use induction by $d=\dim (X)$ in order to show that $\im (q\otimes X)=0$ for any irreducible $X$. If $d=0$ then $\im (q\otimes X)=CH^i(M)$. Since $CH^i(M_{\Omega })=0$ and we work with coefficients in $\QQ $, it follows that $CH^i(M)=0$ too. Suppose we have shown that $\im (q\otimes X)=0$ for $X$ of dimensions $0,1,\dots ,d-1$, and let $X$ be a variety of dimension $d$. Again, $CH ^i(M_K)=0$ because $CH^i(M_{\Omega })=0$, and $\im (q\otimes Z)=0$ by induction hypothesis. Hence $\im (q\otimes X)=0$ as well.

Finally, if $X$ is smooth and projective, we get $\im (q\otimes X)=CH^i(M\otimes M(X))=0$. Then $M=0$ by Manin's principle.
\end{proof}

\section{Representability of algebraic cycles}
\label{s-Fanomotive}

Let $V$ be a smooth projective variety over $k$. The group $A^i_{\ZZ }(V)$ is said to be (weakly) representable, \cite{BM}, if there exists a smooth projective curve $\Gamma$, a cycle class $z$ in $CH^i_{\ZZ }(\Gamma\times V)$, and an algebraic subgroup $G\subset J_{\Gamma }$ in the Jacobian variety $J_{\Gamma }$, such that for any algebraically
closed field $\Omega $ containing $k$ the induced homomorphism
$$
z_*:J_{\Gamma }(\Omega)=A^1_{\ZZ }(\Gamma_{\Omega})\to A^i_{\ZZ }(V_{\Omega})
$$
is surjective, and its kernel is the group $G(\Omega)$. Working with coefficients in $\QQ$ the matter of closeness of $\ker (z_*)$ will be omitted, and representability of $A^i(X)$ means the existence of a surjective homomorphism $z_*$. We call this rational representability.

Recall that for any smooth projective variety $V$, there is a canonical homomorphism
$$
\lambda^i_l: CH^i_{\ZZ}(V)\{l\}\lra H^{2i-1}_{\acute e t}(V,\QQ_l/\ZZ_l(i))\; ,
$$
constructed in \cite{Blo79}, where $CH^i_{\ZZ}(V)\{l\}$ is the $l$-power torsion subgroup in $CH^i_{\ZZ}(V)$. The homomorphisms $\lambda^i_l$ are functorial with respect to the action of correspondences between smooth projective varieties.

The following result is well-know, and we put it here for the convenience of the reader. The argument is a combination of the results of Roitman, Jannsen, Bloch, and Srinivas, see \cite{Jan}.

\begin{lemma}\label{l-represent}
Let $V$ be a smooth projective variety of dimension $d$ over an algebraically closed field $k$. Assume resolution of singularities for varieties of dimension $d-1$ if the field $k$ is of positive characteristic. Suppose that the group $A^d(V)$ is rationally representable. Then the following is true:
\begin{enumerate}
\item
the kernel of the Albanese morphism $A^d_{\ZZ}(V)\to \Alb(V)(k)$ is trivial,
\item
all classes in $H^2(V)(1)$ and $H^{2d-2}(V)(d-1)$ are algebraic,
\item
the integral Chow group $A^2_{\ZZ}(V)$ is representable,
\item
the cokernel of the homomorphism $A^2_{\ZZ }(V)\to CH^2_{\ZZ }(V)_{\rm hom}$ vanishes if $\cha(k)=0$ and is $p$-torsion if $\cha(k)=p>0$. Also, $A^2(V_{\Omega})=CH^2(V_{\Omega})_{hom}$ for any algebraically closed field $\Omega \supset k$,
\item
the homomorphism
$$
\lambda^2_l(V):CH^2_{\ZZ}(V)\{l\}\lra H^3_{\acute e t}(V,\QQ_l/\ZZ_l(2))
$$
is injective and has finite cokernel.
\end{enumerate}
\end{lemma}

\begin{proof}
(1) By definition of rational representability, there exists a smooth projective curve $\Gamma$ and a cycle class $z\in CH^d(\Gamma\times V)$ such that the homomorphism $z_*:A^1(\Gamma)\to A^d(V)$ is surjective. Let $Z$ be a codimension $d$ algebraic cycle on $\Gamma\times V$ representing the class $z$. Suppose $Z'$ is an irreducible component of $Z$ such that the projection $p_\Gamma:Z'\to \Gamma$ is not generically finite. Then the image $p_\Gamma(Z')$ is strictly less than $\Gamma$ and the homomorphism $[Z']_*: A^1(\Gamma)\to A^d(V)$ vanishes. Thus, one can assume that each irreducible component $Z_i$ of the cycle $Z=\sum_i m_iZ_i$ projects generically finite onto $\Gamma$.

It is easy now to compute the action of $z_*$ in terms of symmetric powers of algebraic varieties. For any variety $W$ let $S^nW$ be its $n$-th symmetric power. Let also $a_i$ be the degree of the generically finite morphism $p_\Gamma: Z_i\to \Gamma$, and put $m=\sum_i |m_i|a_i$. Then we have an obvious map
  $$
  \Gamma\times\Gamma \lra S^mV\times S^mV
  $$
by taking a preimage with additional multiplicities $m_i$ of the first point on $\Gamma $ under the surjections $Z_i\to \Gamma $ and then sending everything onto $V$ via the compositions $Z_i\subset \Gamma \times V\to V$, and the same with multiplicities $-m_i$ for the second point on $\Gamma$. For any natural number $n$ we also have a map
  $$
  S^n\Gamma\times S^n\Gamma \lra S^{mn}V\times S^{mn}V\; ,
  $$
and the square
  $$
  \diagram
  S^n\Gamma \times S^n\Gamma \ar[dd]_-{} \ar[rr]^-{} & &
  S^{mn}V\times S^{mn}V \ar[dd]^-{} \\ \\
  A^1_{\ZZ }(\Gamma ) \ar[rr]^-{} & & A^d_{\ZZ }(V)
  \enddiagram
  $$
is commutative, where the vertical maps send $(x_1,\dots ,x_n,y_1,\dots ,y_n)$ to the corresponding sums of differences $x_i-y_i$.
As $\Gamma$ is a curve, there exists a natural number $n$, such that the left vertical mapping is onto. Since the map
$$
z_*:A^1(\Gamma)\to A^d(V)
$$
is surjective by the rational representability of zero-cycles on $V$, we have that the $\QQ$-span of the image of the right vertical map
  $$
  S^{mn}V\times S^{mn}V\lra A^d_{\ZZ }(V)
  $$
is the whole group $A^d(V)$. Moreover, this remains true after we extend scalars to any algebraically closed field $\Omega$ containing $k$, in particular, one can take $\Omega$ to be uncountable. It follows then from \cite[Prop. 1.6]{Jan} that the kernel of the Albanese morphism $A^d_{\ZZ}(V)\to \Alb(V)(k)$ is trivial. \\

(2) The same Proposition 1.6 in \cite{Jan} says that there exists a curve $E\subset V$, such that $CH^d((V-E)_{\Omega})=0$ for any algebraically closed field $\Omega\supset k$. By \cite{BS}, there exists a natural number $N$, a divisor $D\subset V$, and codimension $d$ cycles $R_1,R_2$ on $V\times V$ such that $\supp R_1\subset E\times V$, $\supp R_2\subset V\times D$, and
  $$
  [N\Delta_V]=[R_1]+[R_2]
  $$
in $CH^d(V\times V)$. Following op.cit., consider resolutions of singularities $i':\widetilde{E}\to E$, $j':\widetilde{D}\to D$ and put $i:\widetilde{E}\to V$ and
$j:\widetilde{D}\to V$ be the compositions of $i'$ and $j'$ with the corresponding closed embedings of $E$ and $D$ into $V$. Let also $\widetilde{R}_1$ and $\widetilde{R}_2$ be strict transforms of $R_1$ and $R_2$ on $\widetilde{E}\times V$ and $V\times \widetilde{D}$ respectively.
Using the projection formula for \'etale cohomology and the decomposition $[N\Delta_V]=[R_1]+[R_2]$ it is not hard now to show that any class in $H^2(V)(1)$ is a sum of a class from $(\tilde R_1)_*H^2(\widetilde{E})(1)$ and a class from $j_*H^0(\widetilde{D})$. As $\tilde E$ is a smooth curve, the groups $H^2(\widetilde{E})(1)$ and $H^0(\widetilde{D})$ are algebraic. Then all classes in $H^2(V)(1)$ are algebraic. Since the intersection pairing between algebraic classes in $H^2(V)(1)$ and $H^{2d-2}(V)(d-1)$ is non-degenerate, all classes in $H^{2d-2}(V)(d-1)$ are algebraic too.
\\

(3) It is proved in \cite[Theorem 1(i)]{BS} that under the above assumptions
the integral group $A^2_{\ZZ}(V)$ is representable (even if the characteristic of $k$ is positive) in the sense that there exists an abelian variety $A$ over $k$ and a regular isomorphism of groups $A^2_{\ZZ}(X_{\Omega})\to A(\Omega)$, where $\Omega\supset k$ is a universal domain. Therefore, there exists an algebraic subgroup $G\subset \Pic^0(\widetilde{D})$ such that $G(\Omega)$ is the kernel of surjective homomorphism $j_*:A^1_{\ZZ}(\widetilde{D}_{\Omega})\to A^2_{\ZZ}(V_{\Omega})$. This implies that $A^2_{\ZZ}(V)$ is representable in the sense used in the present paper. \\

(4) The first assertion is contained in \cite[Theorem 1(ii),(iii)]{BS}. It follows immediately that $CH^2(V_{\Omega})_{hom}=A^2(V_{\Omega})$. \\

(5) The injectivity of $\lambda^2_l$ was proved in \cite{MS}, and we need only to show that the cokernel is finite. Since $\widetilde E$ is a smooth curve the group $H^3_{\acute e t}(\widetilde{E},\QQ_l/\ZZ_l(2))$ vanishes. Acting on \'etale cohomology by correspondences $R_1$ and $R_2$ with
$$
[N\Delta_V]=[R_1]+[R_2],
$$
we see that the group $N\cdot H^3_{\acute e t}(V,\QQ_l/\ZZ_l(2))$ coincides with the image of the homomorphism
$$
j_*:H^1_{\acute e t}(\widetilde{D},\QQ_l/\ZZ_l(1))\to H^3_{\acute e t }(V,\QQ_l/\ZZ_l(2))\;.
$$
As it was mentioned above, actions of correspondences commute with the homomorphisms $\lambda _l$. Therefore, as
$$
\lambda^1_l:CH^1_{\ZZ}(\widetilde{D})\{l\}\to H^1_{\acute e t}(\widetilde{D},\QQ_l/\ZZ_l(1))
$$
is bijective by \cite{Blo79}, the image of $\lambda^2_l$ contains the subgroup $N\cdot H^3_{\acute e t}(V,\QQ_l/\ZZ_l(2))$.

Let $r$ be the $l$-adic valuation of $N$. Then we have that
$$
N\cdot H^3_{\acute e t}(V,\QQ_l/\ZZ_l(2))=l^r\cdot H^3_{\acute e t}(V,\QQ_l/\ZZ_l(2))\;.
$$
The exact sequence of \'etale sheaves on $V$
$$
0\to \mu_{l^r}^{\otimes 2}\to \QQ_l/\ZZ_l(2)\stackrel{l^r}\to \QQ_l/\ZZ_l(2)\to 0
$$
yields the exact sequence of cohomology groups
$$
H^3_{\acute e t}(V,\QQ_l/\ZZ_l(2))\stackrel{l^r}\to H^3_{\acute e t}(V,\QQ_l/\ZZ_l(2))\to H^4_{\acute e t}(V,\mu_{l^r}^{\otimes 2})\;.
$$
Since the group $H^4_{\acute e t}(V,\mu_{l^r}^{\otimes 2})$ is finite, we conclude that the cokernel of $\lambda^2_l$ is finite.
\end{proof}

\section{Motivic vision of representability}\label{s-vision}

Now let us have a motivic look at the representability of algebraic cycles in the case of dimension $3$. Let $X$ be a smooth projective threefold over an algebraically closed field $k$. Fix a closed point $x_0$ on $X$ and let
  $$
  \pi_0=[x_0\times X]\; \; \; \; \; \hbox{and}\; \; \; \; \; \pi_6=[X\times x_0]
  $$
be two projectors splitting $\uno $ and $\Le ^3$ from $M(X)$. Let also
  $$
  \pi_1\; \; \; \; \; \hbox{and}\; \; \; \; \; \pi_5
  $$
be the Picard and Albanese projectors respectively, constructed in \cite{Mu1} (see also \cite{Scholl}) and splitting the Picard motive $M^1(X)$ and the Albanese motive $M^5(X)$ from $M(X)$.

Let, furthermore, $\{D'_i\}$, $i=1,\ldots,b$, be a collection of $\QQ $-divisors whose classes give a base in the Neron--Severi group
$NS_{\QQ}^1(X)$. Note that the cycles
$$
D_i=D_i'-\pi_1(D_i')
$$
have the same class in the Neron--Severi group as $D'_i$ and $\pi_1(D_i)=0$. Analogously, let $\{E_i\}$ be rational $1$-cycles giving a base in the dual group $NS^{\QQ}_1(X)$ such that $\pi_1(E_i)=0$ and $\langle D_i\cdot E_j\rangle =\delta_{ij}$. The correspondence
  $$
  \pi _2=\left [\sum _{i=1}^b E_i\times D_i\right ]\in CH^3(X\times X)
  $$
is a projector, and we set
  $$
  \pi _4=\pi_2^{\rm t}\; .
  $$
Then $(X,\pi _2,0)\cong \Le^{\oplus b}$ and $(X,\pi _4,0)\cong(\Le^{ 2})^{\oplus b}$.

All these projectors $\pi _i$ are pairwise orthogonal. This is well-known for $\pi_i$, $i\ne 2,4$, \cite{Scholl}, and the assertion for $\pi_2$ and $\pi_4$ can be proved by the arguments used in \cite[Prop.14.2.3]{KMP}. Therefore, we get a decomposition
$$
M(X)\cong \uno \oplus M^1(X)\oplus \Le ^{\oplus b}\oplus N\oplus (\Le ^{2})^{\oplus b} \oplus M^5(X)\oplus \Le ^3\; ,
$$
where
  $$
  N=(X,\pi ,0)
  $$
and
  $$
  \pi=\Delta _X-\pi _0-\pi _1-\pi _2-\pi _4-\pi _5-\pi _6\; .
  $$

Notice that for any $i\in \{ 0,1,5,6\} $ the projector $\pi _i$ acts identically on $H^i(X)$ and trivially on $H^j(X)$ when $j\ne i$. The motivic Poincar\'e duality for $M(X)$ (see~\cite{Scholl}) induces an isomorphism $N\cong N^{\vee}\otimes\Le^3$.

Assume now that $A^3(X)$ is rationally representable. Then for any index $i=2,4$ the projector $\pi _i$ acts trivially on $H^j(X)$ if $j\ne i$ and identically on $H^i(X)$, because $NS_{\QQ_l}^1(X)=H^2(X)(1)$ and $NS^{\QQ_l}_1(X)=H^4(X)(2)$ by Lemma \ref{l-represent} (2). Besides, all the projectors $\pi_i$, $i\in \{ 0,1,2,5,6\} $, act trivially on $CH^2(X)$ (see \cite[Theorems 1,\,2]{Mu1}) and \cite[Theorem 4.4 (iii)]{Scholl}, and the action of $\pi_4$ on $CH^2(X)$ is given by the formula:
$$
\pi_4(a)=\sum_{i=1}^b\langle a \cdot D_i\rangle E_i
$$
for any $a\in CH^2(X)$, so that the kernel of the action of $\pi _4$ coincides with the group $CH^2(X)_{\rm num}$.

The Chow group $CH^i(N)$ is, by definition, the image of the action of the correspondence $\pi $ on $CH^i(X)$. Since that correspondence is the orthogonal complement in the diagonal to the sum $\pi _0+\pi _1+\pi _2+\pi _4+\pi _5+\pi _6$ it follows that $CH^i(N)$ is, at the same time, the kernel of the action of this sum.

It is not hard to see that $CH^i(N)=0$ when $i\neq 2$. As all $\pi_i$, $i\in \{ 0,1,2,5,6\} $, act trivially on $CH^2(X)$, we have that $CH^2(N)$ is nothing but the kernel of the action of $\pi _4$, which is the group $CH^2(X)_{\rm num}$. On the other hand, $CH^1(X)_{\rm num}=CH^1(X)_{\rm hom}$ for divisors. As the group $H^2(X)(1)$ is algebraic by Lemma \ref{l-represent}(2) and the pairing between $H^2(X)(1)$ and $H^4(X)(2)$ is non-degenerate, it follows that $CH^2(X)_{\rm num}$ coincides with $CH^2(X)_{\rm hom}$, and $CH^2(X)_{\rm hom}=A^2(X)$ by Lemma~\ref{l-represent}(4). Thus,
  $$
  CH^2(N)=A^2(X)\; .
  $$

Clearly,
  $$
  H^3(N)=H^3(X)\; \; \; \; \hbox{and}\; \; \; \; H^j(N)=0\; \; \hbox{for}\; j\ne 3\; .
  $$

It all means that, provided $A^3(X)$ is representable, $N$ is the middle motive in $M(X)$ controlling algebraic cycles in codimension $2$ on $X$. Our main aim is to show that $N$ is finite-dimensional bringing finite-dimensionality for the entire motive $M(X)$, see Theorem~\ref{t-zerocycles} and Corollary~\ref{findim} below.

By Lemma \ref{l-represent} (3), the integral Chow group $A^2_{\ZZ}(X)$ is representable. Let $\Gamma$ be a smooth projective curve, $z\in CH^2_{\ZZ}(\Gamma \times X)$ be a correspondence and let $G$ be an algebraic subgroup in $J_{\Gamma }$ coming from the definition of the integral representability of $A^2_{\ZZ}(X)$. The cycle class $z$ can be considered also as a morphism
$$
z:M_{\ZZ }(\Gamma)\otimes \Le \lra M_{\ZZ}(X)\; .
$$
The quotient $J=J_{\Gamma }/G$ is an abelian variety, and the corresponding projection
$$
\alpha :J_{\Gamma }\lra  J
$$
induces a morphism
$$
M(\alpha ):M(J)\lra M(J_{\Gamma })\; .
$$
Fixing a point on the curve $\Gamma$ we obtain an embedding
$$
i_\Gamma:\Gamma\hra J_{\Gamma }\; ,
$$
which gives a morphism
$$
M(i_\Gamma):M(J_{\Gamma })\lra M(\Gamma)\; .
$$
Composing all of these three morphisms, twisted by $\Le $ when appropriate, we obtain a morphism
$$
w=z\circ(M(\alpha\circ i_\Gamma)\otimes id_{\Le}):M_{\ZZ}(J)\otimes\Le\to M_{\ZZ}(X).
$$

\begin{lemma}\label{l-Chowisogeny}
The integral correspondence $w: M_{\ZZ }(J)\otimes \Le \to M_{\ZZ }(X)$
induces a homomorphism of integral Chow groups
$$
w_*: A^1_{\ZZ }(J)\lra A^2_{\ZZ }(X)\; ,
$$
which is an isogeny, i.e. is surjective with finite kernel.
\end{lemma}

\begin{proof}
Recall that for any abelian variety $A$ there exists an ample line bundle $\mathcal L$ on the dual variety $A^{\vee }$ defining an isogeny $\varphi _{\mathcal L}:A^{\vee }\to A$, see \cite{Mum}. For $J_{\Gamma }$ there exists an ample line bundle $\Theta $ on $J_{\Gamma }^{\vee }$, so-called theta-bundle, such that $\varphi _{\Theta }:J_{\Gamma }^{\vee }\to J_{\Gamma }$ is an isomorphism of abelian varieties, which coincides with the homomorphism $i_\Gamma^*:A^1_{\ZZ }(J_{\Gamma })\to A^1_{\ZZ}(\Gamma)$ via the natural identifications of $J_{\Gamma }^{\vee }$ with $A^1_{\ZZ}(J_{\Gamma })$ and $J_{\Gamma }$ with $A^1_{\ZZ}(\Gamma)$. The surjective morphism $\alpha :J_{\Gamma }\to J$ induces an injective dual morphism $\alpha ^{\vee }:J^{\vee }\hra J_{\Gamma }^{\vee }$, and the composition $J^{\vee}\stackrel{\alpha^{\vee}}\lra J_{\Gamma }^{\vee}\stackrel{\varphi_{\Theta}}\lra J_{\Gamma }\stackrel{\alpha}\lra J$
coincides with the isogeny $\varphi _{\mathcal L}$, where $\mathcal L$ is a restriction of $\Theta $ on $J^{\vee}$, loc. cit. The identification $J^{\vee}=A^1_{\ZZ}(J)$ shows that the composition of the group homomorphisms
$$
A^1_{\ZZ}(J)\stackrel{\alpha^*}\lra A^1_{\ZZ}(J_{\Gamma })\stackrel{i_{\Gamma}^*}\lra
A^1_{\ZZ}({\Gamma})=J_{\Gamma}\stackrel{\alpha }\lra J
$$
is an isogeny in the above sense. Since $w_*=z_*\circ i^*_{\Gamma }\circ \alpha ^*$ and $z_*: A^1_{\ZZ }({\Gamma })\to A^2_{\ZZ }(X)$ factors through $J$, such that $J\cong A^2_{\ZZ }(X)$, we get the needed assertion.
\end{proof}

For any abelian group $A$ let $A\{ l\} $ be the subgroup of $l$-primary torsion in $A$, i.e. the union of the kernels of all multiplication by $l^i$ homomorphisms $A\stackrel{l^i}{\to }A$.

\begin{corollary}
\label{extra}
The homomorphism
$$
w_*:CH^1_{\ZZ}(J)\{l\}\to CH^2_{\ZZ}(X)\{l\}
$$
has finite kernel and cokernel.
\end{corollary}

\begin{proof}
By Lemma \ref{l-Chowisogeny} the homomorphism $A^1_{\ZZ }(J)\lra A^2_{\ZZ }(X)$ is surjective with finite kernel. It follows that the homomorphism $A^1_{\ZZ }(J)\{ l\} \lra A^2_{\ZZ }(X)\{ l\} $ has finite kernel and is surjective, too. Consider the following commutative square:
  $$
  \diagram
  A^1_{\ZZ }(J)\{ l\}  \ar[dd]_-{} \ar[rr]^-{w_*} & & A^2_{\ZZ }(X)\{ l\}
  \ar[dd]^-{} \\ \\
  CH^1_{\ZZ }(J)\{ l\} \ar[rr]^-{w_*} & & CH^2_{\ZZ }(X)\{ l\}
  \enddiagram
  $$
The left vertical arrow is an isomorphism because $J$ is an abelian variety. The right vertical arrow decomposes as
  $$
  A^2_{\ZZ }(X)\{ l\}\to CH^2_{\ZZ }(X)_{\rm hom}\{ l\}\to
  CH^2_{\ZZ }(X)\{ l\} \; .
  $$
The cokernel of the injective homomorphism $A^2_{\ZZ }(X)\to CH^2_{\ZZ }(X)_{\rm hom}$ is $p$-torsion by Lemma \ref{l-represent}(4). Since $l\neq p$ the first homomorphism in the above composition is an isomorphism $A^2_{\ZZ }(X)\{ l\} \cong CH^2_{\ZZ }(X)_{\rm hom}\{ l\} $. The cokernel of the second injective homomorphism in the composition is contained in the group $H^4_{\acute e t}(X,\ZZ _l(2))\{l\}$, so is finite.
\end{proof}

We need some more $l$-adic tools. Recall that for any smooth projective variety $V$ of dimension $d$ over $k$ we have a nondegenerate bilinear pairing
  \begin{equation}
  \label{zetmodn}
  H^j_{\acute e t}(V,\ZZ /n(i))\times H^{2d-j}_{\acute e t}(V,\ZZ /n(d-i))\lra
  \ZZ /n\; .
  \end{equation}
When $n=l^s$ we can pass to the limit of the system
  $$
  \dots \to \ZZ /l^{s+1}\lra \ZZ /l^s\to \dots \; .
  $$
getting a pairing
  \begin{equation}
  \label{zetel}
  H^j_{\acute e t}(V,\ZZ _l(i))\times H^{2d-j}_{\acute e t}(V,\ZZ _l(d-i))\lra
  \ZZ _l\; .
  \end{equation}

The point is that we can also pass to the colimit of the system
  $$
  \dots \to \ZZ /l^s\ZZ \stackrel{\times l}{\lra }\ZZ /l^{s+1}\ZZ \to \dots \; ,
  $$
for the second variable in (\ref{zetmodn}) having a pairing
  \begin{equation}
  \label{raska}
  H^j_{\acute e t}(V,\ZZ _l(i))\times H^{2d-j}_{\acute e t}(V,\QQ _l/\ZZ _l(d-i))\lra
  \QQ _l/\ZZ _l\; ,
  \end{equation}
which is also nondegenerate, see \cite[Theorem 1.11]{Raskind}.

A straightforward tensoring of (\ref{zetel}) with $\QQ _l/\ZZ _l$ gives a pairing
  $$
  H^j_{\acute e t}(V,\ZZ _l(i))\times (H^{2d-j}_{\acute e t}(V,\ZZ _l(d-i))\otimes \QQ _l/\ZZ _l)\lra \QQ _l/\ZZ _l\; ,
  $$
which induces a homomorphism
  $$
  H^{2d-j}_{\acute e t}(V,\ZZ _l(d-i))\otimes \QQ _l/\ZZ _l\lra \Hom (H^j_{\acute e t}(V,\ZZ _l(i)),\QQ _l/\ZZ _l)\; .
  $$
The nondegenerate pairing (\ref{raska}) gives an isomorphism
  $$
  \Hom (H^j_{\acute e t}(V,\ZZ _l(i)),\QQ _l/\ZZ _l)\cong H^{2d-j}_{\acute e t}(V,\QQ _l/\ZZ _l(d-i))\; ,
  $$
loc.cit. Thus, one has a homomorphism
  $$
  \varrho :H^{2d-j}_{\acute e t}(V,\ZZ _l(d-i))\otimes \QQ _l/\ZZ _l\lra H^{2d-j}_{\acute e t}(V,\QQ _l/\ZZ _l(d-i))\; .
  $$

Since the pairing $(\ref{zetel})$ has finite kernels on both variables, the homomorphism $\varrho $ has finite kernel and cokernel.

\begin{lemma}
\label{zayac}
The homomorphism $w_*:H^1(J)(-1)\to H^3(X)$ is bijective.
\end{lemma}

\begin{proof}
Consider the commutative diagram
$$
\diagram
CH^1_{\ZZ}(J)\{l\} \ar[dd]^-{\lambda^1_l} \ar[rr]^-{w_*}& &  CH^2_{\ZZ}(X)\{l\}
\ar[dd]^-{\lambda^2_l} \\ \\
H^1_{\acute e t}(J,\QQ_l/\ZZ_l(1)) \ar[rr]^-{w_*} & &
H^3_{\acute e t}(X,\QQ_l/\ZZ_l(2))
\enddiagram
$$
The homomorphism $\lambda_l^1$ is bijective by \cite{Blo79}. The homomorphism $\lambda^2_l$ is injective and has finite cokernel by Lemma~\ref{l-represent}(5). The top horizontal homomorphism has finite kernel and cokernel by Corollary~\ref{extra}. Therefore, the bottom horizontal homomorphism has finite kernel and cokernel too.

The diagram
$$
\diagram
H^1_{\acute e t}(J,\ZZ_l(1))\otimes_{\ZZ_l}\QQ_l/\ZZ_l \ar[dd]^-{\varrho }  \ar[rr]^-{w_*}& &  H^3_{\acute e t}(X,\ZZ_l(2))\otimes_{\ZZ_l}\QQ_l/\ZZ_l
\ar[dd]^-{\varrho } \\ \\
H^1_{\acute e t}(J,\QQ_l/\ZZ_l(1)) \ar[rr]^-{w_*} & &
H^3_{\acute e t}(X,\QQ_l/\ZZ_l(2))
\enddiagram
$$
commutes and, as we have seen just now, the vertical homomorphisms $\varrho $ have finite kernels and cokernels. Since the bottom horizontal homomorphism has finite kernel and cokernel, so is the top horizontal homomorphism. Since \'etale cohomology groups with $\ZZ _l$-coefficients are finitely generated $\ZZ _l$-modules, the homomorphism
   $$
   w_*:H^1_{\acute e t}(J,\ZZ_l(1))\lra H^3_{\acute e t}(X,\ZZ_l(2))
   $$
also has finite kernel and cokernel (it is, in fact, injective but we do not need that). Therefore, after tensoring with $\QQ _l$, the homomorphism $w_*:H^1(J)(-1)\to H^3(X)$ is bijective.
\end{proof}

\begin{remark}
{\rm If $\cha(k)=0$ then one can use relations between abelian varieties and polarizable weight one Hodge structures getting a more direct proof of Lemma~\ref{zayac} that will not use the homomorphism $\lambda^2_l$.}
\end{remark}

The decomposition of rational motives of abelian varieties, \cite{DM}, \cite{Kuen},
$$
M(J)\cong \bigoplus_{i=0}^{2\dim(J)}\wedge^i M^1(J)\; ,
$$
yields an injective morphism
$$
i_J:M^1(J)\lra M(J)\; .
$$
We put
$$
M=M^1(J)\otimes \Le \; ,
$$
and let
$$
f:M\lra M(X)
$$
be the composition $f=\pi \circ w\circ i_J$. Then $f$ can be also considered as a morphism into the middle motive $N$, i.e.
  $$
  f:M\lra N\; .
  $$

\begin{corollary}
\label{l-inejectivecohom}
The correspondence $f:M\to N$ induces an bijective homomorphism
$$
f_*: H^3(M)=H^1(M^1(J))(-1)\stackrel{\cong}\lra H^3(N)\; .
$$
\end{corollary}

\begin{proof}
The homomorphism
  $$
  w_*:H^1(J)(-1)\to H^3(X)
  $$
is bijective by Lemma \ref{zayac}. Since
  $$
  (i_J)_*:H^1(M^1(J))\lra H^1(M(J))=H^1(J)
  $$
  is
an isomorphism, and $\pi $ induces an isomorphism between $H^3(X)$ and $H^3(N)$, the composition $f=\pi\circ w\circ i_J$ induces an isomorphism
$$
f_*: H^1(M^1(J))(-1)\stackrel{(i_J)_*}\lra H^1(J)(-1)\stackrel{w_*}\lra H^3(X)
\stackrel{\pi _*}\lra H^3(N)\; .
$$

Let now $\pi _1^J$ be the Picard projector for the variety $J$, so that $M^1(J)=(J,\pi _1^J,0)$. Then,
  $$
  H^1(M^1(J))(-1)=H^1(J,\pi ^J_1,0)(-1)=H^3(J,\pi ^J_1,-1)=H^3(M)\; .
  $$
This completes the proof of the lemma.
\end{proof}

\section{Motives of threefolds with representable $A^3(X)$}

In this section we will prove a theorem giving a precise motivic criterium of weak representability of algebraic cycles on a threefold over an arbitrary algebraically closed field $k$.

\begin{theorem}
\label{t-zerocycles}
Let $X$ be a smooth projective threefold over an algebraically closed field $k$. The group $A^3(X)$ is rationally representable if and only if the motive $M(X)$ has the following Chow-K\"unneth decomposition:

$$
M(X)\cong \uno\oplus M^1(X)\oplus\Le^{\oplus b}\oplus
(M^1(J)\otimes\Le )\oplus(\Le ^2)^{\oplus b}\oplus M^5(X)\oplus \Le ^3\; ,
$$

\medskip

\noindent where $M^1(X)$ and $M^5(X)$ are the Picard and Albanese motives respectively, $b=b^2(X)=b^4(X)$ is the Betti number, and $J$ is a certain abelian variety over $k$, isogenous to the intermediate Jacobian $J^2(X)$ if $k=\CC $.
\end{theorem}

\begin{proof}
Notice that if $M(X)$ has a decomposition of the above type then $A^3(X)$ is representable. Indeed, the only Chow group of a Picard motive is concentrated in codimension one, \cite[Theorem 1]{Mu1}. Therefore, from the above decomposition we get
  $$
  A^3(X)\cong A^3(M^5(X))\; .
  $$
The hard Lefschetz theorem for Chow motives gives an isomorphism
  $$
  M^1(X)\otimes \Le ^2\stackrel{\cong }{\lra }M^5(X)\; ,
  $$
\cite[Theorem 4.4 (ii)]{Scholl}. Then
  $$
  A^3(M^5(X))\cong A^1(M^1(X))\; ,
  $$
and the last thing is representable since it is a Chow group of codimension one algebraic cycles.

Conversely, assume $A^3(X)$ is representable. We use notations and constructions from Section~\ref{s-vision}. Our strategy is to split the motive $M$ from $N$ through the morphism $f$ and to show that the rest has no Chow groups, so vanishes by Lemma \ref{l-gersten}. Then the middle motive $N$ will be isomorphic to the motive $M=M^1(J)\otimes \Le $, and combining all the gadgets developed above we will obtain the desired Chow-K\"unneth decomposition.

First notice that in any pseudo-abelian category $\bcA $ if $A\stackrel{f}{\to}B$ and $B\stackrel{g}{\to}A$ are two morphisms, such that the composition $gf$ is an automorphism $h$ of the object $A$, then $fh^{-1}g$ is an idempotent on $B$.

In our case, we want to construct a morphism $g:N\to M$ such that the composition $g\circ f$ would be an automorphism of the motive $M$. That can be done with the aid of Corollary \ref{l-inejectivecohom}. Indeed, let $n=\dim (J)$ and let $M^{2n-1}(J)$ be the Albanese motive for the abelian variety $J$. Since the Albanese projector is nothing but the transposition of the Picard projector, we have that
  $$
  M^{2n-1}(J)=M^1(J)^{\vee }\otimes \Le ^n\; .
  $$
By the hard Lefschetz theorem for Chow motives,
  $$
  M^1(J)\otimes \Le ^{n-1}\stackrel{\cong }{\lra }M^{2n-1}(J)\; ,
  $$
\cite[Theorem 4.4 (ii)]{Scholl}. Then we have the following isomorphism:
  $$
  M^{\vee }\otimes \Le ^3=M^1(J)^{\vee }\otimes \Le ^2=
  M^{2n-1}(J)\otimes \Le ^{2-n}\cong M^1(J)\otimes \Le =M\; .
  $$
Define $g$ to be a composition
  $$
  \diagram
  g:N=N^{\vee }\otimes \Le ^3 \ar[rr]^-{f^{\vee }\otimes \id _{\Le ^3}} & &
  M^{\vee }\otimes \Le ^3 \ar[r]^-{\cong } & M\; .
  \enddiagram
  $$
As $f$ induces an isomorphism on the third cohomology groups by Corollary \ref{l-inejectivecohom}, the morphism $g$ also induces an isomorphism on $H^3$,
  $$
  g_*:H^3(N)\stackrel{\cong }{\lra }H^3(M)\; .
  $$
We need to show that the composition
  $$
  h=gf
  $$
is an automorphism of the motive $M$. There is a canonical isomorphism
$$
\End (M^1(J))=\End (J)\otimes_{\ZZ }\QQ \; ,
$$
where $\End (J)$ is the ring of regular endomorphisms of the abelian variety $J$, see \cite[Proposition 4.5]{Scholl}. Each such an endomorphism induces an endomorphism of the first cohomology of $J$, so we have a homomorphism
$$
\End (J)\otimes_{\ZZ }\QQ \lra \End (H^1(J))\; ,
$$
Let $\tilde h$ be an endomorphism of $J$ such that its rational multiple is induced by $h$. Since $h$ induces an automorphism on the group $H^3(M)=H^1(J)(-1)$ the endomorphism $\tilde h$ is an isogeny. Let $\tilde h^{-1}$ be its quasi-inverse isogeny, so that the composition $\tilde h^{-1}\tilde h$ is a multiplication by a natural number, and the same for $\tilde h\tilde h^{-1}$. Tensoring with $\QQ $ we see that $\tilde h$ is invertible in the algebra $\End (J)\otimes_{\ZZ }\QQ $. Hence, $h$ is an automorphism of the motive $M$.

Then immediately we have a splitting
  $$
  N=M\otimes N'
  $$
in $\CHM(k)$ defined by the projector
  $$
  \pi _3=fh^{-1}g
  $$
and its orthogonal complement
  $$
  \pi _3'=\pi -\pi _3\; .
  $$
Then
  $$
  \Delta _X=\pi _3'+\sum _{i=0}^6\pi _i
  $$
is a decomposition of the diagonal into eight pair-wise orthogonal idempotents in the associative ring $\End(M(X))$. Respectively, for each index $j=0,1,2,3$ we have a decomposition
  $$
  CH^j(X)=CH^j(X,\pi _3',0)\oplus \left(\oplus _{i=0}^6CH^j(X,\pi _i,0)\right)\; .
  $$
It is well-known that
  $$
  \oplus _{i=0}^6CH^j(X,\pi _i,0)=CH^j(X)
  $$
for all $j\in \{ 0,1,3\} $, see \cite[Theorem 4.4 (iii)]{Scholl}, and the same holds over any algebraically closed extension $\Omega \supset k$. As to codimension $2$, we know already that
$$
w_*: A^1(J)\to A^2(X)
$$
is an isomorphism by Lemma~\ref{l-Chowisogeny}, which remains valid after a scalar extension $\Omega\supset k$. Since
$$
(i_J)_*:CH^1(M^1(J)_{\Omega})\stackrel{\cong }\lra A^1(J_{\Omega})
$$
is an isomorphism and $\pi $ acts identically on $A^2(X_{\Omega})$, we see that $f$ defines an isomorphism of rational Chow groups
$$
f_*: CH^2(M_{\Omega})=CH^1(M^1(J_{\Omega}))\stackrel{\cong }\lra A^2(X_{\Omega})\; .
$$
Since $\pi_3=fh^{-1}g$,
$$
CH^2(M_{\Omega})=CH^2((X,\pi_3,0)_{\Omega})=A^2(X_{\Omega})\; .
$$
Moreover,
$$
CH^2((X,\pi_4,0)_{\Omega})=NS^{\QQ}_1(X_{\Omega})\; .
$$
Then,
$$
CH^2((X,\pi _3+\pi _4,0)_{\Omega})=CH^2(X_{\Omega})\; ,
$$
so that
  $$
  CH^j(N'_{\Omega })=0
  $$
for any $j$. Then $N'=0$ by Lemma \ref{l-gersten}.

By Corollary~\ref{l-inejectivecohom}, the correspondence $w:M(J)\otimes\Le\to M(X)$ gives an isomorphism $w_*:H^1(J)(-1)\to H^3(X)$. If $k=\CC$ this gives a corresponding isomorphism of rational Hodge structures, so that the intermediate Jacobian $J^3(X)$ is algebraic. Since the integral correspondence $w:M_{\ZZ}(J)\otimes\Le\to M_{\ZZ}(X)$ acts on integral Hodge structures, it acts also on intermediate Jacobians, i.e. there exists a morphism of abelian varieties $w_*:J\to J^2(X)$. Since $w_*$ induces an isomorphism on first rational cohomology groups of these abelian varieties, it is an isogeny. This completes the proof.
\end{proof}

Recall that a Chow motive is said to be abelian if it is an object in the full tensor pseudo-abelian subcategory generated by motives of curves in $\CHM (k)$. By Kimura's results, all abelian motives are finite-dimensional, see \cite{Kimura}.

\begin{corollary}
\label{findim}
Let $X$ be a smooth projective variety of dimension $3$ over an algebraically closed field $k$, such that the group $A^3(X)$ is rationally representable. Then the motive $M(X)$ is finite-dimensional.
\end{corollary}

\begin{proof}
All the components in the above decomposition of $M(X)$ are abelian, so finite-dimensional.
\end{proof}

\section{Threefolds fibred by surfaces with algebraic $H^2$}
\label{s-fibred}

By the results of Bloch, Murre and Koll\'ar the group $A^3(X)$ is representable if $X$ is a Fano threefold over an algebraically closed field of characteristic zero. In this section we will prove representability (with coefficients in $\QQ $) of $A^3(X)$ for threefolds fibred over a curve with
representable algebraic cycles in the generic fibre and satisfying a certain additional condition. We will use desingularization of algebraic surfaces, which is known to be always valid, so that we proceed to work over an algebraically closed field $k$ of arbitrary characteristic.

We need to introduce some more terminology and prove a certain localization lemma. First we should recall a localized version of weak representability. Let $V$ be a smooth (not necessary projective) variety over $k$, and let $A$ be a quotient of $A^i(V)$, i.e., we are given a surjective morphism of $\QQ$-vector spaces
  $$
  A^i(V)\lra A\; .
  $$
Then $A$ is said to be representable if there exists a smooth projective curve $\Gamma $ and a cycle class $z$ in $CH^i(\Gamma \times V)$, such that the composition
  $$
  A^1(\Gamma)\stackrel{z_*}\lra A^i(V)\to A
  $$
is surjective.

Notice that, as we use Chow groups with coefficients in $\QQ $, we omit the matter of closeness of the kernel of the homomorphism
$z_*$.

In particular, given a relative Chow motive $M$ over a smooth base $S$, one can speak about representability of the group $A^i(M)$ because it is a quotient of $A^i(V)$ for a certain smooth $S$-variety $V$. It is easy to show that representability of $A^i(M)$ does not depend on the choice of a variety $V$ and a projector on it.

Let $V$ be an open subset in a smooth projective variety $W$. The restriction homomorphism $A^i(W)\to A^i(V)$ is surjective, and representability of $A^i(V)$ as an identical quotient of itself is equivalent to representability of $A^i(V)$ as a quotient of $A^i(W)$.

If $A^i(V)$ splits into a finite direct sum of subspaces, each of which is representable, then the whole group $A^i(V)$ is representable too. This was established in \cite{Bloch} (see somewhat more detailed expansion of that idea in \cite{Gul}). Therefore, in proving representability, we may cover $A^i(V)$ by zero-cycles on a finite collection of smooth projective algebraic curves, and then patch altogether into a representability of $A^i(V)$ by a unique curve $\Gamma $.

We will also need a method of tracking information on algebraic cycles from the generic fibre of a family of varieties to a fibre over a closed point. This can be done in a few ways, actually. Possibly the simplest and most visual one is through so-called ``spread+limit" construction in smooth families.

Let $R$ be a discrete valuation ring with the maximal ideal $\gm $, $K$ be the fraction field and $l=R/\gm $ be the residue field. Let $\bar K$ and $\bar l$ be algebraic closures of the fields $K$ and $l$ respectively. Let also $S=\Spec (R)$. For any smooth morphism of schemes $V\to S$, we have a
a specialization homomorphism
  $$
  \sigma :CH ^i(V_K)\lra CH ^i(V_l)
  $$
which commutes with intersection products, pull-backs and push-forwards of algebraic cycles, see \cite[20.3]{Ful}.

In geometrical situation $\sigma $ can be interpreted as follows. Let $C$ be an algebraic curve of an algebraically closed field $k$, and let $R$ be the local ring at a closed point on it. If $X\to C$ is a smooth family over $C$, then $\sigma $ gives a homomorphism from the Chow group $CH^i(X_{\eta })$ of the geometrical generic fibre of the family $X\to C$ to the Chow group $CH^i(X_{p})$ of its fibre over a closed point $p$ on $C$. Let $Z$ be an algebraic cycle on $X_{\eta }$. Let $Z'$ be the Zariski closure of a spread of $Z$ over some Zariski open subset in $C$ containing the point $p$. Then
  $$
  \sigma ([Z])=[Z'_p]\; .
  $$

We need to know two things about geometrical specialization homomorphisms. First of all, specializations $\sigma $ commute with intersection products, pull-backs and push-forwards of algebraic cycles. It follows immediately that, given a smooth family $X\to C$ over a curve $C$, finite-dimensionality of the motive of the generic fibre $X_{\eta }$ (or $X_{\bar \eta }$) implies finite-dimensionality of the motive of the fibre $X_p$ over a closed point $p$ on $C$. Indeed, motivic finite-dimensionality of $X_\eta$ is equivalent to the decomposition of the diagonal $\Delta _{X_{\eta}}$ into two orthogonal projectors $\alpha_{\eta}$ and $\beta_{\eta}$, such that a wedge power of $\alpha_{\eta}$ and a symmetric power of $\beta_{\eta}$ are rationally trivial on $X_{\eta }^N$ for some $N$. Applying the specialization $\sigma$ we obtain that the same property for the diagonal holds true in the special fibre $X_p$. Certainly, this can be also proved in a more sophisticated way by means of Ayoub's motivic vanishing cycle functor, see \cite{Ayoub}.

Another thing is that specializations commute with cycle maps, \cite[20.3]{Ful}. More precisely, if $R$ is Henselian, say complete, then passing to colimits over finite extensions of $R$ we also have a specialization homomorphism over algebraically closed fields:
  $$
  \sigma :CH ^i(V_{\bar K})\lra CH ^i(V_{\bar l})\; .
  $$
Moreover, one has a commutative diagram
  $$
  \diagram
  CH^i(V_{\bar K}) \ar[dd]_-{\sigma } \ar[rr]^-{cl} & & H^{2i}(V_{\bar K})(i) \ar[dd]^-{\sigma \, ,\; \cong } \\ \\
  CH^i(V_{\bar k}) \ar[rr]^-{} & & H^{2i}(V_{\bar k})(i)
  \enddiagram
  $$

\begin{lemma}
\label{newbrick}
Let $X\to C$ be a smooth morphism of relative dimension two onto an algebraic curve $C$ over $k$. Let $p$ be a closed point on $C$. Then, if $A^2(X_{\bar \eta })$ is representable, the group $A^2(X_p)$ is also representable.
\end{lemma}

\begin{proof}
Since zero-cycles on the surface $X_{\bar \eta }$ are representable, it follows that its diagonal is balanced by curves, \cite{BS}, so the motive $M(X_{\bar \eta })$ is finite-dimensional, \cite{GP2}, and the group $H^2(X_{\bar \eta})(1)$ is algebraic, \cite{TheBible}. Using a henselisation of the local ring at $p$ and the above commutative square we see that the algebraicity of $H^2(X_{\bar\eta })(1)$ implies algebraicity of $H^2(X_p)(1)$. Moreover, the finite-dimensionality of $M(X_{\bar \eta })$ implies finite-dimensionality of $M(X_p)$ through the specialization $\sigma $. Then $A^2(X_p)$ is representable by \cite{GP2}.
\end{proof}

\begin{lemma}
\label{local}
Let $X\to C$ be a projective dominant morphism from a smooth irreducible threefold $X$ over $k$ onto a smooth irreducible curve $C$ over $k$. Let $U$ be a non-empty Zariski open subset in $C$, and let $Y=X\times _CU$ be the preimage of $U$ under the morphism $X\to C$. Then $A^2(X)$ is representable if and only if $A^2(Y)$ is representable. If, moreover, zero-cycles in the geometric generic fibre $X_{\bar \eta }$ and in desingularizations of all the degenerate fibers of the morphism $X\to C$ are representable, then the same assertion holds in codimension three: $A^3(X)$ is representable if and only if $A^3(Y)$ is representable.
\end{lemma}

\begin{proof}
Let $S$ be a finite collection of points on the curve $C$, such that $U$ is a complement of $S$ in $C$. Let $i$ be either $2$, or $3$. By the localization sequence,
  $$
  \oplus _{p\in S}CH^{i-1}(X_p)\lra CH^i(X)\lra CH^i(Y)\to 0\; ,
  $$
the $\QQ $-vector space $CH^i(X)$ splits into two parts,
  $$
  CH^i(X)=CH^i(Y)\oplus I\; ,
  $$
where $I$ is the image of the left hand side homomorphism. As it was mentioned above, the right hand side surjection induces a surjective homomorphism on algebraically trivial algebraic cycles,
  $$
  A^i(X)\lra A^i(Y)\; ,
  $$
and we get the splitting
  $$
  A^i(X)=A^i(Y)\oplus A\; ,
  $$
where $A=I\cap A^i(X)$ (the intersection is being taken inside $CH^i(X)$) after we fix a section $A^i(Y)\to A^i(X)$ and identify $A^i(Y)$ with its image in $A^i(X)$. It follows immediately that representability of $A^i(X)$ implies representability of $A^i(Y)$ and, in order to derive representability of $A^i(X)$ from representability of $A^i(Y)$ all we need is to show representability of $A$, which depends on representability of $A^{i-1}(X_p)$ for $p\in S$.

For codimension $i=2$ the argument is simple. For any $p\in S$ let $\tilde X_p$ be a resolution of singularities of the surface $X_p$ ($\tilde X_p=X_p$ if $X_p$ is smooth). We have that $A^i(\tilde X_p)=\oplus A^i(B)$, where the sum is taken over all irreducible components $B$ of the surface $\tilde X_p$. The groups $A^1(B)$ are representable, hence the group $A^1(\tilde X_p)$ is representable. The N\'eron-Severi group $NS^1_{\QQ }(B)=CH^1(B)/A^1(B)$ is finite-dimensional. In addition, each morphism $\tilde X_p\to X_p$ induces a surjective push-forward homomorphism on Chow groups with rational coefficients. It follows that the kernel $A$ of the localization homomorphism $A^2(X)\to A^2(Y)$ splits into two $\QQ $-vector subspaces $A_1$ and $A_2$, where $A_1$ representable and the second subspace $A_2$ is finite-dimensional, so representable too.

In codimension $3$ the situation is slightly more subtle. If $X_p$ is a smooth fiber, then $A^2(X_p)$ is representable by Lemma \ref{newbrick}. If $X_p$ is a degenerate fibre $X\to C$, then we need representability of $A^2(\tilde X_p)$ for a desingularization $\tilde X_p$ of the degenerate fibre $X_p$, but this is what we have by the assumption of the lemma. In both cases, $A^2(\tilde X_p)$ is representable, and the quotient $CH^2(\tilde X_p)/A^2(\tilde X_p)$ is finite dimensional. It follows that $A$ is representable.
\end{proof}

\begin{remark}
\label{r-local}
{\rm The same argument as in the proof of Lemma~\ref{local} implies the following statement. Let $M=(V/S,p,0)$ be a relative Chow motive over a smooth base $S$. Recall that it means, in particular, that the morphism $V\to S$ is smooth and projective, and $p$ is a relative projector on $V$. Then, for any non-empty open subset $U\subset C$, representability of $A^2(M)$ is equivalent to representability of $A^2(M_U)$, where $M_U=(V_U/U,p_U,0)$, $V_U=V\times_C U$, and $p_U$ is the restriction of $p$ to $V_U$.}
\end{remark}

\bigskip

Let again $X\to C$ be a projective dominant morphism from a smooth irreducible threefold $X$ over $k$ onto a smooth irreducible curve $C$.
Denote by $J$ the relative Albanese variety of the smooth part of $X\to C$. There is a smooth morphism $J\to C^0$, where $C^0\subset C$ is the open subset over which the morphism $X\to C$ is smooth. Let $M^1(J/C^0)$ be the relative Picard projector for the relative abelian variety $J$ over $C^0$, see \cite{Kuen}.

\begin{theorem}
\label{main2}
Let $X\to C$ be a projective dominant morphism from a smooth irreducible
threefold $X$ over $k$ onto a smooth irreducible curve $C$ over $k$,
and let $J$ and $M^1(J/C^0)$ be as above. Suppose that
$A^2(X_{\bar \eta })$ is representable for the geometric generic
fibre $X_{\bar \eta }$ and that $A^2(M^1(J/C^0))$ is representable.
Then $A^2(X)$ is representable. If, in addition, zero-cycles
in desingularizations of all the degenerate fibers of the
morphism $X\to C$ are representable, then $A^3(X)$ is representable.
\end{theorem}

\begin{proof}
Let $U$ be a Zariski open subset in $C$ and let $Y=X\times_C U$ be its preimage. By Lemma~\ref{local}, it is sufficient to show representability of $A^i(Y)$ for $i=2$ and $3$.

Murre's projectors for a surface can be actually constructed over a non-algebraically closed field, \cite{Mu1}, \cite{Scholl}. Consider the Murre
decomposition of the diagonal for the generic fibre $X_{\eta}$,
  $$
  \Delta _{X_{\eta}}=\pi _0+\pi _1+\pi _2+\pi _3+\pi _4,
  $$
which gives the corresponding decomposition for the motive
  $$
  M(X_{\eta })=\oplus _{j=0}^4M^j(X_{\eta })\; ,
  $$
see \cite{Mu1}, \cite{Scholl}. Recall that $M^0(X_{\eta })=\uno $, $M^1(X_{\eta })\cong M^1(J_{\eta})$ by \cite[Proposition 4.5]{Scholl}, $M^3(X_{\eta })\cong M^1(J_{\eta })\otimes \Le $, $M^4(X_{\eta })=\Le^{2}$ and the motive $M^2(X_{\eta })$ controls the Albanese kernel of $X_{\eta }$.

Recall that, by \cite{GP2}, representability of $A^2(X_{\bar\eta})$ is equivalent to the fact that all classes in $H^2(X_{\bar\eta})(1)$ are algebraic and the motive $M(X_{\bar \eta})$ is finite-dimensional. Working with rational coefficients, finite-dimensionality of $M(X_{\bar \eta})$ is equivalent to finite-dimensionality of $M(X_{\eta})$. Moreover, the group $H^2(X_{\bar\eta})(1)$ is algebraic. Therefore, we have that
  $$
  M^2(X_{\bar\eta })=\Le^{\oplus b^2}\; ,
  $$
where $b^2=\dim H^2(X_{\bar\eta })$ is the second Betti number of $X_{\bar\eta }$, see \cite{GP2}.

Let now $\Pi _i$ be a spread of $\pi _i$ as an algebraic cycle, $i=0,\ldots ,4$. Cutting out fibers of the morphism $X\to C$ we get that all $\Pi _i$ are idempotents too. Notice that by Lemma~\ref{local} representability of $A^i(X)$ is equivalent to representability of $A^i(Y)$ for $i=2,3$. So, we have relative motives
  $$
  M^i(Y/U)=(Y/U,\Pi _i,0)
  $$
in the category $\CHM (U)$ over $U$. Respectively, one has
the decomposition
  $$
  M(Y/U)=\bigoplus_{j=0}^4M^j(Y/U)
  $$
in $\CHM (U)$. Since
  $$
  \CHM (\eta )=\colim _V\CHM (V)\; ,
  $$
where the colimit is taken over all non-empty open subsets $V$ in
$C$, we see that making $U$ smaller we obtain isomorphisms of relative
motives
  $$
  M^0(Y/U)\cong \uno \; ,
  $$
  $$
  M^1(Y/U)\cong M^1(J_U/U)\; ,
  $$
  $$
  M^3(Y/U)\cong M^1(J_U/U)\otimes \Le \; ,
  $$
  $$
  M^4(Y/U)\cong \Le^2
  $$
and
  $$
  M^2(Y'/U')\cong \Le^{\oplus b_2}\; ,
  $$
where $U'\to U$ is a certain finite cover, $Y'=Y\times_UU'$, and $J_U=J\times_{C^0}U$.

\medskip

The above decomposition of the relative Chow motive $M(Y/U)$ implies the corresponding decomposition of the Chow groups $CH^i(Y)$ for $i=2$ and $3$:

$$
CH^i(Y)=
$$
$$
=CH^i(U)\oplus CH^i(M^1(J_U/U))\oplus I\oplus CH^{i-1}(M^1(J_U/U))\oplus CH^{i-2}(U)\; ,
$$

\bigskip

\noindent where $I$ is a direct summand in the group $CH^{i-1}(U')^{\oplus b_2}$. By Corollary 3.2 in \cite{DM} the group $CH^3(M^1(J_U/U))$ consists of all elements $\alpha\in CH^3(J_U)$ such that for any natural number $n$ we have $n^*\alpha=\alpha$, where $n:J_U\to J_U$ is the multiplication by $n$. Thus, $CH^3(M^1(J_U/U))=CH^3_5(J_U,\QQ)$ in notations from op.cit. Therefore by Theorem~2.19 in op.cit. the group $CH^3(M^1(J_U/U))$ vanishes.

\medskip

Clearly, the analogous decomposition holds for the subgroup $A^i(Y)$ in $CH^i(Y)$ generated by cycles algebraically equivalent to zero. Therefore, we have that
  $$
  A^2(Y)=A^2(M^1(J_U/U))\oplus A\oplus A^1(M^1(J_U/U))\; ,
  $$
and
  $$
  A^3(Y)=A^2(M^1(J_U/U))\oplus A^1(U)\; ,
  $$

\medskip

\noindent where $A$ is a direct summand in $A^1(U')^{\oplus b_2}$.

The groups
  $$
  A^1(U')^{\oplus b_2}\; ,\; \; A\; ,\; \; A^1(M^1(J_U/U))\; \; \; \; \hbox{and}\; \; \; A^1(U)
  $$
are always representable. Therefore, the representability of $A^i(Y)$ for $i=2,3$ is equivalent to representability of the group $A^2(M^1(J_U/U))$. By Remark~\ref{r-local}, the latter is equivalent to representability of $A^2(M^1(J/C^0))$. This completes the proof.
\end{proof}

\section{Applications}

The purpose of this section is to apply the above general results in order to prove motivic finite-dimensionality for certain classes of threefolds $X$ over $k$. From now on we will assume that $k$ is an algebraically closed field of characteristic zero.

The first thing to say is that motivic finite-dimensionality is a birational invariant in dimension three:

\begin{lemma}
Let $X$ and $Y$ be two smooth projective threefolds over $k$. Assume that $X$ is  birationally equivalent to $Y$. Then $M(X)$ is finite-dimensional if and only if $M(Y)$ is finite-dimensional.
\end{lemma}

\begin{proof}
By Hironaka's theorem, there is a blow-up with smooth centers $\widetilde{X}\to X$ such that the induced rational morphism $\widetilde{X}\to Y$ is regular. As motives of smooth curves are finite-dimensional, Manin's motivic formula for blow ups, \cite[2.7]{Scholl}, and the standard properties of finite-dimensional objects in tensor categories show that finite-dimensionality of $X$ implies finite-dimensionality of $Y$. By symmetry between $X$ and $Y$, this proves the needed result.
\end{proof}

The following theorem gives unconditional results on motivic finite-dimensionality for Fano threefolds and threefolds fibred by Del Pezzo or Enriques surfaces over a curve. Notice that any Fano threefold can be fibred by $K3$ surfaces, and $K3$ surfaces cover Enriques' ones.

\begin{theorem}
\label{findimconcrete}
Let $X$ be a smooth projective threefold over $k$, which belongs to one of the following two types of varieties:

\begin{itemize}

\item
smooth projective Fano threefolds;

\item
up to birational equivalence $X$ is fibred by Del Pezzo or Enriques surfaces over a curve.

\end{itemize}

\noindent Then the group $A^3(X)$ is representable, the motive $M(X)$ is finite-dimen\-sional, and it has the Chow-K\"unneth decomposition as in Theorem \ref{t-zerocycles}.
\end{theorem}

\begin{proof}
Let first $X$ be a Fano threefold over $k$. It is well know that the groups $A^*(X)$ are all representable on $X$, see \cite{Mu3}. Then $M(X)$ has the desired decomposition by Theorem \ref{t-zerocycles}, and so it is finite-dimensional by Corollary \ref{findim}.

Assume now that, up to a birational equivalence, $X$ is fibred over a smooth irreducible curve $C$, $X\to C$, with generic fibre being Enriques' surface. The motive of an Enriques surface is known to be finite-dimensional, see \cite{GP1}. Moreover, an Enriques surface is regular. Therefore in notations of Theorem~\ref{main2}, the relative Albanese variety $J$ vanishes. Thus, by Theorem \ref{main2}, representability of $A^3(X)$ depends on representability of zero cycles in desingularizations of degenerate fibers of the morphism $X\to C$. Let $p$ be a point on $C$, such that the corresponding fibre $X_p$ is degenerate. By Mumford's semistable reduction theorem, without loss of generality one can assume that $X_p$ is a normal crossing divisor on $X$. By Kulikov's classification of semi-stable degeneration of Enriques' surfaces \cite{Kulikov}, we have that each component in $X_p$ belongs to the following three types: smooth Enriques, elliptic ruled and rational. In all of these cas
 es the geometrical genus of each irreducible component in $X_p$ is equal to zero and the motive is finite-dimensional. Hence, the corresponding Chow group of zero-cycles is representable. From this we see that the Chow group of zero cycles of the desingularization of $X_p$ is representable too. Then $A^3(X)$ is representable by Theorem \ref{main2}. Now, again, $M(X)$ has the desired decomposition by Theorem \ref{t-zerocycles}, and is finite-dimensional by Corollary \ref{findim}. The arguments for $X$ fibred by Del Pezzo surfaces are analogous.
\end{proof}

One can also state a conditional result which still have concrete geometrical meaning. Below $q$ stands for the irregularity of a surface.

\begin{theorem}
\label{findimconcrete2}
Let $X$ be a smooth projective threefold over $k$ up to a birational equivalence fibred over a curve with the generic fibre being a surface with $p_g=0$ and $q=1$. Let, furthermore, $E$ be a desingularization of a projective closure of a spread of the elliptic curve $\Alb ^0(X_{\eta })$, where $X_{\eta}$ is the generic fibre. Then, again, $A^3(X)$ is representable, the motive $M(X)$ is finite-dimensional and it has the Chow-K\"unneth decomposition as in Theorem \ref{t-zerocycles}, provided $p_g(E)=0$.
\end{theorem}

\begin{proof}
Indeed, in terms of Theorem \ref{main2} the surface $E$ is just a desingularization of a projective closure of the relative Albanese variety $J$. As $p_g(E)=0$ and $E$ is an elliptic surface, the motive $M(E)$ is finite-dimensional, \cite{GP1}. It is equivalent to say that $A^2(E)$ is representable. It follows that $A^2(J)$ is representable and, in particular, in $A^2(M^1(J/C^0))$ is representable too. From results of Morrison, \cite{Morr}, we claim that all irreducible components of a semistable degeneration of a hyperelliptic surface has representable zero-cycles. Then, by Theorem \ref{main2}, the Chow group $A^3(X)$ is representable. Hence, the motive $M(X)$ can be decomposed by Theorem \ref{t-zerocycles}, and it is finite-dimensional by Corollary \ref{findim}.
\end{proof}

In conclusion one can say also that the methods above give us a possibility to prove motivic finite-dimensionality for some threefolds fibred by surfaces of general type with $p_g=0$, provided we know representability of zero-cycles in semi-stable degenerations of the generic fibre. And the same is true for fiberings by ruled surfaces. However, such results are too much conditional, so we do not state them here.

\bigskip

{\it Acknowledgements.} We wish to thank Aleksandr Pukhlikov for several useful discussions and Joseph Ayoub for explaining the meaning of the motivic vanishing functor in our approach. We are also grateful to Dmitri Orlov and Ivan Panin for pointing out a mistake in the previous version of the paper. The second named author is grateful to the Institute for Advanced Study in Princeton, where the insight of this work had been partially originated a few years ago.

\bigskip

\bibliographystyle{amsplain}

\end{document}